\documentclass[11pt]{article}
\usepackage[T1]{fontenc}
\usepackage[utf8]{inputenc}
\usepackage{lmodern}
\usepackage{microtype}
\usepackage{amsmath,amssymb,amsthm,mathtools}
\usepackage[margin=2.4cm]{geometry}
\usepackage{booktabs}
\usepackage{parskip}           % no paragraph indentation
\usepackage[backref=page]{hyperref}
\hypersetup{colorlinks=true,linkcolor=blue,urlcolor=blue}

\hypersetup{hidelinks}

\newtheorem{theorem}{Theorem}[section]
\numberwithin{equation}{section}
\newtheorem{proposition}[theorem]{Proposition}

\newtheorem{corollary}[theorem]{Corollary}
\theoremstyle{definition}

\newtheorem{remark}[theorem]{Remark}
\newtheorem{example}[theorem]{Example}

\usepackage{enumerate}          % needed for \begin{enumerate}[(1)]
\newcommand{\selabel}[1]{\label{se:#1}}
\newcommand{\seref}[1]{Section~\ref{se:#1}}

\def\thlabel#1{\label{th:#1}}
\def\thref#1{Theorem~\ref{th:#1}}

\newcommand{\prlabel}[1]{\label{pr:#1}}
\newcommand{\prref}[1]{Proposition~\ref{pr:#1}}
\newcommand{\colabel}[1]{\label{co:#1}}
\newcommand{\coref}[1]{Corollary~\ref{co:#1}}
\newcommand{\relabel}[1]{\label{re:#1}}
\newcommand{\reref}[1]{Remark~\ref{re:#1}}
\newcommand{\exlabel}[1]{\label{ex:#1}}

\newcommand{\eqlabel}[1]{\label{eq:#1}}
\newcommand{\equref}[1]{(\ref{eq:#1})}

\newcommand{\bZ}{\mathbb{Z}}
\newcommand{\bQ}{\mathbb{Q}}

\title{Centralizers of maps as automorphism groups of Yang--Baxter solutions}

\author{G. Militaru}
\date{}

\begin{document}
\maketitle

\begin{abstract}
For a partition $X = Y \sqcup Z$ of a set $X$ and an arbitrary map $g \colon Y \to {\rm Hom}_{\rm Set} (Z, \, Z)$, we construct a solution $R_g \colon X \times X \to X \times X$ of the braid equation. We characterize intrinsically the solutions which are isomorphic to some $R_g$ and we give a necessary and sufficient condition for two solutions $R_g$ and $R_{g'}$ to be isomorphic. Consequently, the number of isomorphism classes of the solutions $R_g$ on a set with $n$ elements is $1 + \sum_{k = 1}^{n-1} \widetilde{o} (n - k, \, k)$, where $\widetilde{o} (m, \, k)$ is the number of orbits of the symmetric group $S_m$, acting by conjugation, on the set of multisets with $k$ elements of maps $\{1, \ldots, m\} \to \{1, \ldots, m\}$. This construction is then used to obtain a solution of the braid equation whose automorphism group is exactly the centralizer ${\rm Aut} (f_1, \ldots, f_k)$ of a given family $f_1, \ldots, f_k \colon X \to X$ of maps, $k \geq 1$. For $k = 2$, combined with a result of Agore, Chirvasitu and the author, this gives another proof of the fact that any group is isomorphic to the automorphism group of a solution of the braid equation (equivalently, of the set-theoretic Yang--Baxter equation).
\end{abstract}

\medskip


\vspace{.3cm}


%%%%%%%%%%%%%%%%%%%%%%%%%%%%%%%%%%%%%%%%%%%%%%%%%%%%%%%%%%%%%%%%%%%%%%%%%%%%%
\section*{Introduction}

The Yang--Baxter equation, in its linear form, is a basic equation of mathematical physics: its solutions are essential for building integrable models in statistical mechanics \cite{Baxter, Yang}. However, the classification of these solutions is an extremely difficult problem: only partial results are known even in dimension $3$ (see \cite{HMR, LR} and the references therein).
For this reason, Drinfel'd \cite{Dr} made a strategic move: he proposed to study first the set-theoretic version of the equation, in which the vector space is replaced by a set and the linear map by a map $R \colon X \times X \to X \times X$. Looking back, and in view of \cite{MS}, the most natural framework for this equation is probably an arbitrary braided monoidal category.

Even at the level of sets, starting with the foundational papers \cite{etingof, LuYZ}, the study of this equation has generated a great deal of good mathematics in various directions: combinatorics and graph theory, knot theory, the theory of (semi)groups and Hopf--Galois extensions, to name only a few. For more details we refer to the recent monograph \cite{BES} and the references therein.

In this paper we continue, from a different angle, the study of the problems recently considered in \cite{MilAut, MilBIP}. First, for a partition $X = Y \sqcup Z$ of a set $X$ and a family $g = (g_y)_{y\in Y}$ of maps $g_y \colon Z \to Z$, we construct in \thref{Tg} a new class of idempotent solutions $R_g \colon X \times X \to X \times X$ of the braid equation; we solve the isomorphism problem for these solutions and we compute their automorphism groups. \prref{recogTg} describes intrinsically the solutions $(X, \, R)$ which are isomorphic to some $(X, \, R_g)$. As an application, the number of isomorphism classes of solutions of the form $R_g$ on a set with $n$ elements is
\begin{equation*}
1 + \sum_{k = 1}^{n-1} \widetilde{o} (n - k, \, k),
\end{equation*}
where $\widetilde{o} (m, \, k)$ is the number of orbits of the symmetric group $S_m$, acting by conjugation, on the set of multisets with $k$ elements of
maps $\{1, \ldots, m\} \to \{1, \ldots, m\}$ (\coref{countTg}). As we show in \reref{countTg}, the numbers $\widetilde{o} (m, \, k)$ can be computed as the coefficients of an explicit generating function, using the Cauchy--Frobenius Lemma and P\'olya's enumeration theorem.

Next, let $f = (f_1, \ldots, f_k)$ be a finite family of maps $f_i \colon X \to X$, $k \geq 1$, and let
\begin{equation*}
{\rm Aut} (f_1, \ldots, f_k) := \bigl\{ \sigma \in {\rm Sym} (X) \, \bigl| \, \sigma \circ f_i = f_i \circ \sigma, \,\, \text{for all} \,\, i = 1, \ldots, k \bigr\}
\end{equation*}
be its \emph{centralizer} \cite[Definition 1.1]{acm1}, which is a subgroup of ${\rm Sym} (X)$. For $k = 1$, it is known that ${\rm Aut} (f_1)$ is the automorphism group of a solution of the braid equation on the set $X$ itself \cite[Theorem 2.19]{acm1}. Using the solutions $R_g$, we show that, for any $k \geq 1$, ${\rm Aut} (f_1, \ldots, f_k)$ is the automorphism group of an idempotent solution of the braid equation on a set with $|X| + 2k + 3$ elements (\thref{centr}). Since any group is the centralizer of a family of two maps \cite[Theorem 4.2(1)]{acm1}, this gives another proof of the fact that any group is the automorphism group of a solution of the braid equation \cite{MilAut}, which can moreover be chosen idempotent (\coref{centr}). Results of this type, going back to Frucht's theorem for graphs, have been obtained in many categories; see, for instance, \cite{BF, Jones} and the references therein.

Finally, the problem of finding, for a given finite group $G$, the smallest set on which $G$ is the automorphism group of a solution of the braid equation remains open (\reref{smallest}).

%%%%%%%%%%%%%%%%%%%%%%%%%%%%%%%%%%%%%%%%%%%%%%%%%%%%%%%%%%%%%%%%%%%%%%%%%%%%%
%%%%%%%%%%%%%%%%%%%%%%%%%%%%%%%%%%%%%%%%%%%%%%%%%%%%%%%%%%%%%%%%%%%%%%%%%%%%%
\section{Preliminaries}\selabel{prel}

Given a set $X$, we denote by $\tau_X \colon X\times X \to X \times X$ the flip map, $\tau_X (x, \, y) = (y, \, x)$; ${\rm Sym} (X)$ denotes the symmetric group on $X$, i.e. the group of all bijections $X \to X$. We recall that a \emph{multiset} is a set in which elements may be repeated, e.g. $\{x, \, x, \, x, \, y, \, \dots \}$. More rigorously, a multiset is an ordered pair $(U, \, \mu)$, where $U$ is a usual set and $\mu \colon U \to \mathbb{N}^{*}$ is a function which gives, for each $u \in U$, the \emph{multiplicity} $\mu (u)$ of $u$ in the multiset; if $U$ is finite, the number of elements of the multiset is $\sum_{u \in U} \mu (u)$. The symmetric group $S_m$ acts by conjugation on the multisets of maps $\{1, \ldots, m\} \to \{1, \ldots, m\}$: for $\beta \in S_m$, the multiset $\beta \cdot (U, \, \mu)$ has the underlying set $\{\beta \circ f \circ \beta^{-1} \, | \, f \in U\}$, and $\beta \circ f \circ \beta^{-1}$ has the multiplicity $\mu (f)$, for all $f \in U$.

For a map $R\colon X \times X \to X\times X$ we use the usual notation
\begin{eqnarray*}
  && R^{12} ,\,\, R^{23},\,\, R^{13} : X \times X \times X \to X \times X \times X\\
  && R^{12} := R \times {\rm Id}_X, \quad  R^{23} := {\rm Id}_X \times R, \quad R^{13}: = ({\rm Id}_X \times \tau_X) \circ R^{12} \circ ({\rm Id}_X \times \tau_X)
\end{eqnarray*}
A map $R: X\times X \to X\times X$ is a solution of the \emph{Yang--Baxter (YB) equation} (resp. \emph{braid equation}) if
\begin{equation}\label{eq:qybe}
R^{12} R^{13} R^{23} = R^{23}R^{13}R^{12}, \qquad  ({\rm resp.} \,\, R^{12}R^{23}R^{12} = R^{23}R^{12}R^{23})
\end{equation}
as maps $X\times X \times X \to X\times X \times X$ under the usual composition. We recall that $R$ is a solution of the braid equation if and only if $\tau_X \circ R$ is a YB-solution. Any map $R \colon X \times X \to X \times X$ can be written as
\begin{equation}\eqlabel{notatieR}
R(x, \, y) = (x \cdot y, \, x \ast y),
\end{equation}
where $(X, \, \cdot, \, \ast)$ is a bi-magma, i.e. a set with two binary operations; this gives a bijection between maps $R \colon X \times X \to X \times X$ and bi-magma structures on $X$. A morphism of solutions $\sigma \colon (X, \, R) \to (X', \, R')$ is a map $\sigma \colon X \to X'$ such that $(\sigma \times \sigma) \circ R = R' \circ (\sigma \times \sigma)$, i.e. a morphism of the associated bi-magmas: $\sigma (x \cdot y) = \sigma(x) \cdot \sigma(y)$ and $\sigma (x \ast y) = \sigma(x) \ast \sigma(y)$. For $x \in X$ we denote by $l_x$, $r_x \colon X \to X$ the \emph{left} and \emph{right translations} of the bi-magma $(X, \, \cdot, \, \ast)$, i.e. $l_x (y) := x \cdot y$ and $r_x (y) := y \ast x$, for all $y \in X$; thus $R (x, \, y) = \bigl( l_x (y), \, r_y (x) \bigr)$.

Throughout we denote by $\mathcal{YB}$ the category of all YB-solutions, i.e. objects of $\mathcal{YB}$ are pairs $(X, \, R)$, where $X$ is a set and $R$ a YB-solution on $X$. We denote by $\mathcal{BR}$ the category of solutions of the braid equation, i.e. objects of $\mathcal{BR}$ are pairs $(X, \, R)$, where $X$ is a set and $R$ is a solution of the braid equation on $X$. The functor
\begin{equation*}\label{iso_br_yb}
\mathcal{YB} \to \mathcal{BR}, \qquad (X, \, R) \mapsto (X, \, \tau_X \circ R)
\end{equation*}
is obviously an isomorphism of categories with inverse given by $(X, \, R) \mapsto (X, \, \tau_X \circ R)$.

The automorphism group ${\rm Aut} (X, \, R)$ of $(X, \, R)$ is the group of all bijections $\sigma \colon X \to X$ with $(\sigma \times \sigma) \circ R = R \circ (\sigma \times \sigma)$, i.e. the automorphism group of the bi-magma $(X, \, \cdot, \, \ast)$. Since $(\sigma \times \sigma) \circ \tau_X = \tau_X \circ (\sigma \times \sigma)$ for any map $\sigma$, we have
\begin{equation}\eqlabel{auttau}
{\rm Aut} (X, \, R) = {\rm Aut} (X, \, \tau_X \circ R).
\end{equation}
All the solutions constructed in this paper are solutions of the braid equation, but they can be equivalently described as YB-solutions, as $(X, \, R) \mapsto (X, \, \tau_X \circ R)$ is an isomorphism between the categories $\mathcal{BR}$ and $\mathcal{YB}$. Therefore, all the results below can be stated for the YB-solutions $r := \tau_X \circ R$.

%%%%%%%%%%%%%%%%%%%%%%%%%%%%%%%%%%%%%%%%%%%%%%%%%%%%%%%%%%%%%%%%%%%%%%%%%%%%%

\section{A new family of solutions of the braid equation: recognition and counting}\selabel{Tg}

In this section we construct, for a partition $X = Y \sqcup Z$ of a set $X$ and a family $g = (g_y)_{y\in Y}$ of maps $g_y \colon Z \to Z$, an idempotent solution $R_g$ of the braid equation on $X$; then we characterize the solutions isomorphic to some $R_g$ and we count them on a finite set. The construction is close to the one of the bipartite solutions of \cite[Theorem 2.2]{MilBIP}, but the maps $g_y$ by which the elements of $Y$ act on $Z$ are no longer required to be bijective.

\begin{theorem}\thlabel{Tg}
Let $Y$ and $Z$ be two disjoint sets, $g = (g_y)_{y\in Y}$ a family of maps $g_y \colon Z \to Z$ and $X := Y \sqcup Z$. Let $R_g \colon X \times X \to X \times X$ be the map given by
\begin{equation}\eqlabel{Tg}
R_g (x, \, x') := \begin{cases} \bigl( g_y (z), \, y \bigr), & \text{if} \,\,\, x = y \in Y, \,\, x' = z \in Z,\\
(x, \, x'), & \text{otherwise}.
\end{cases}
\end{equation}
Then:
\begin{enumerate}[(a)]
\item $R_g$ is an idempotent solution of the braid equation.
\item Let $(Y', \, Z', \, g')$ be another triple as above and $X' := Y' \sqcup Z'$. If $Y$, $Z$, $Y'$, $Z'$ are non-empty, then $(X, \, R_g) \cong (X', \, R_{g'})$ if and only if there exist two bijections $\alpha \colon Y \to Y'$ and $\beta \colon Z \to Z'$ such that
\begin{equation}\eqlabel{isoTg}
\beta \circ g_y = g'_{\alpha(y)} \circ \beta, \qquad \text{for all} \,\, y \in Y.
\end{equation}
\item If $Y \neq \emptyset \neq Z$, then the set
\begin{equation*}
H_g := \Bigl\{ (\alpha, \, \beta) \in {\rm Sym} (Y) \times {\rm Sym} (Z) \, \bigl| \, \beta \circ g_y = g_{\alpha(y)} \circ \beta, \,\, \text{for all} \,\, y\in Y \Bigr\}
\end{equation*}
is a subgroup of ${\rm Sym} (Y) \times {\rm Sym} (Z)$ and the map $H_g \to {\rm Aut} (X, \, R_g)$, $(\alpha, \, \beta) \mapsto \alpha \sqcup \beta$, is an isomorphism of groups.
\end{enumerate}
\end{theorem}

\begin{proof}
(a) Since $R_g (y, \, z) = (g_y (z), \, y) \in Z \times Y$ and $R_g$ is the identity on $Z \times Y$, we have $R_g^2 = R_g$. Write $R := R_g$. If a triple $(x_1, \, x_2, \, x_3) \in X^3$ contains no element of $Y$ followed immediately by an element of $Z$, i.e. if it is of one of the forms $(y, \, y', \, y'')$, $(z, \, z', \, z'')$, $(z, \, y, \, y')$ or $(z, \, z', \, y)$, then all the maps $R^{12}$ and $R^{23}$ fix it, and both sides of the braid equation fix it. In the remaining four cases, for $y$, $y'\in Y$ and $z$, $z' \in Z$, the two sides of the braid equation are computed as follows:
\begin{align*}
(y, \, y', \, z) &\overset{R^{12}}{\longmapsto} (y, \, y', \, z) \overset{R^{23}}{\longmapsto} \bigl( y, \, g_{y'} (z), \, y' \bigr) \overset{R^{12}}{\longmapsto} \bigl( g_y g_{y'} (z), \, y, \, y' \bigr),\\
(y, \, y', \, z) &\overset{R^{23}}{\longmapsto} \bigl( y, \, g_{y'} (z), \, y' \bigr) \overset{R^{12}}{\longmapsto} \bigl( g_y g_{y'} (z), \, y, \, y' \bigr) \overset{R^{23}}{\longmapsto} \bigl( g_y g_{y'} (z), \, y, \, y' \bigr);\\[6pt]
(y, \, z, \, y') &\overset{R^{12}}{\longmapsto} \bigl( g_y (z), \, y, \, y' \bigr) \overset{R^{23}}{\longmapsto} \bigl( g_y (z), \, y, \, y' \bigr) \overset{R^{12}}{\longmapsto} \bigl( g_y (z), \, y, \, y' \bigr),\\
(y, \, z, \, y') &\overset{R^{23}}{\longmapsto} (y, \, z, \, y') \overset{R^{12}}{\longmapsto} \bigl( g_y (z), \, y, \, y' \bigr) \overset{R^{23}}{\longmapsto} \bigl( g_y (z), \, y, \, y' \bigr);\\[6pt]
(y, \, z, \, z') &\overset{R^{12}}{\longmapsto} \bigl( g_y (z), \, y, \, z' \bigr) \overset{R^{23}}{\longmapsto} \bigl( g_y (z), \, g_y (z'), \, y \bigr) \overset{R^{12}}{\longmapsto} \bigl( g_y (z), \, g_y (z'), \, y \bigr),\\
(y, \, z, \, z') &\overset{R^{23}}{\longmapsto} (y, \, z, \, z') \overset{R^{12}}{\longmapsto} \bigl( g_y (z), \, y, \, z' \bigr) \overset{R^{23}}{\longmapsto} \bigl( g_y (z), \, g_y (z'), \, y \bigr);\\[6pt]
(z, \, y, \, z') &\overset{R^{12}}{\longmapsto} (z, \, y, \, z') \overset{R^{23}}{\longmapsto} \bigl( z, \, g_y (z'), \, y \bigr) \overset{R^{12}}{\longmapsto} \bigl( z, \, g_y (z'), \, y \bigr),\\
(z, \, y, \, z') &\overset{R^{23}}{\longmapsto} \bigl( z, \, g_y (z'), \, y \bigr) \overset{R^{12}}{\longmapsto} \bigl( z, \, g_y (z'), \, y \bigr) \overset{R^{23}}{\longmapsto} \bigl( z, \, g_y (z'), \, y \bigr).
\end{align*}
In all cases the two sides coincide, so $R = R_g$ is a solution of the braid equation.

(b) We show that a bijection $\sigma \colon X \to X'$ is an isomorphism of solutions $(X, \, R_g) \to (X', \, R_{g'})$ if and only if $\sigma = \alpha \sqcup \beta$, for two bijections $\alpha \colon Y \to Y'$ and $\beta \colon Z \to Z'$ satisfying \equref{isoTg}, where $\alpha \sqcup \beta$ is the bijection which is equal to $\alpha$ on $Y$ and to $\beta$ on $Z$.

Let $\alpha$ and $\beta$ be two bijections satisfying \equref{isoTg} and $\sigma := \alpha \sqcup \beta$. Since $\sigma (Y) = Y'$ and $\sigma (Z) = Z'$, the maps $(\sigma \times \sigma) \circ R_g$ and $R_{g'} \circ (\sigma \times \sigma)$ are both equal to $\sigma \times \sigma$ outside $Y \times Z$. For $y \in Y$ and $z \in Z$ we have
\begin{equation*}
(\sigma \times \sigma) \bigl( R_g (y, \, z) \bigr) = \bigl( \beta (g_y (z)), \, \alpha(y) \bigr) = \bigl( g'_{\alpha(y)} (\beta(z)), \, \alpha(y) \bigr) = R_{g'} \bigl( \sigma(y), \, \sigma(z) \bigr).
\end{equation*}
Hence $\sigma$ is an isomorphism of solutions.

Conversely, let $\sigma \colon X \to X'$ be an isomorphism of solutions. By \equref{Tg}, the pairs which are not fixed by $R_g$ are exactly the elements of $Y \times Z$:
\begin{equation*}
R_g (x, \, x') \neq (x, \, x') \,\, \Longleftrightarrow \,\, (x, \, x') \in Y \times Z,
\end{equation*}
since the first component of $R_g (y, \, z) = (g_y (z), \, y)$ belongs to $Z$, for all $y\in Y$ and $z \in Z$; similarly, the pairs which are not fixed by $R_{g'}$ are exactly the elements of $Y' \times Z'$. As $R_{g'} \circ (\sigma \times \sigma) = (\sigma \times \sigma) \circ R_g$ and $\sigma \times \sigma$ is injective, a pair $(x, \, x')$ is fixed by $R_g$ if and only if $(\sigma(x), \, \sigma(x'))$ is fixed by $R_{g'}$. Thus, $(x, \, x') \in Y \times Z$ if and only if $(\sigma(x), \, \sigma(x')) \in Y' \times Z'$. Now, if $y \in Y$ and $z \in Z$, then $(\sigma (y), \, \sigma (z)) \in Y' \times Z'$; hence $\sigma (Y) \subseteq Y'$ and $\sigma (Z) \subseteq Z'$. As $\sigma$ is a bijection from $X = Y \sqcup Z$ to $X' = Y' \sqcup Z'$, we obtain $\sigma (Y) = Y'$ and $\sigma (Z) = Z'$, i.e. $\sigma = \alpha \sqcup \beta$, for two bijections $\alpha \colon Y \to Y'$ and $\beta \colon Z \to Z'$. Finally, for $y \in Y$ and $z \in Z$ we have
\begin{equation*}
\bigl( \beta (g_y (z)), \, \alpha(y) \bigr) = (\sigma \times \sigma) \bigl( R_g (y, \, z) \bigr) = R_{g'} \bigl( \sigma (y), \, \sigma (z) \bigr) = \bigl( g'_{\alpha(y)} (\beta(z)), \, \alpha(y) \bigr),
\end{equation*}
i.e. $\alpha$ and $\beta$ satisfy \equref{isoTg}.

(c) It is easy to see that $H_g$ is a subgroup of ${\rm Sym} (Y) \times {\rm Sym} (Z)$ and that $(\alpha, \, \beta) \mapsto \alpha \sqcup \beta$ is an injective morphism of groups. By the proof of (b), applied for $(Y', \, Z', \, g') = (Y, \, Z, \, g)$, its image is ${\rm Aut} (X, \, R_g)$.
\end{proof}

The literature contains a large number of constructions of solutions of the braid equation (equivalently, of the Yang--Baxter equation) from various categories or classes of algebraic and combinatorial objects. As we pointed out in \cite[Introduction, Problem B]{ACM2}, for each such construction it is important to know exactly which solutions are isomorphic to the ones in the given family, i.e. to find an intrinsic characterization of these solutions, expressed only in terms of the map $R$. Unfortunately, very few results of this type are known, although we consider this problem essential. For the solutions $R_g$ given by \equref{Tg}, such a characterization is given by the following:

\begin{proposition}\prlabel{recogTg}
Let $X$ be a set and $R \colon X \times X \to X \times X$ a map. Then $(X, \, R)$ is isomorphic to $(Y \sqcup Z, \, R_g)$, for some disjoint sets $Y$, $Z$ and some family $g = (g_y)_{y \in Y}$ of maps
$g_y \colon Z \to Z$, if and only if there exists a subset $A \subseteq X$ such that, for all $x$, $x' \in X$:
\begin{equation}\eqlabel{recogTg}
R (x, \, x') \neq (x, \, x') \,\, \Longleftrightarrow \,\, x \in A \,\, \text{and} \,\, x' \notin A, \quad \text{and in this case} \quad R (x, \, x') \in (X \setminus A) \times \{x\}.
\end{equation}
Moreover, if $R \neq {\rm Id}_{X \times X}$, then the subset $A$ is unique.
\end{proposition}

\begin{proof}
Assume first that $A \subseteq X$ satisfies \equref{recogTg} and let $Z_A := X \setminus A$. For $a \in A$ and $z \in Z_A$ we have $R (a, \, z) = \bigl( l_a (z), \, r_z (a) \bigr) \in Z_A \times \{a\}$, i.e. $l_a (z) \in Z_A$ and $r_z (a) = a$. Hence, for any $a \in A$, the left translation $l_a$ induces a map $h_a \colon Z_A \to Z_A$, $h_a (z) := l_a (z)$, and, by \equref{recogTg} and \equref{Tg}, $R = R_h$, where $h = (h_a)_{a \in A}$; in particular, $R$ is an idempotent solution of the braid equation, by \thref{Tg}(a). If $R \neq {\rm Id}_{X \times X}$, then $A$ is the set of all $x \in X$ such that $R (x, \, x') \neq (x, \, x')$ for some $x' \in X$, hence it is unique.

Conversely, let $\sigma \colon X \to Y \sqcup Z$ be a bijection such that $(\sigma \times \sigma) \circ R = R_g \circ (\sigma \times \sigma)$, and let $A := \sigma^{-1} (Y)$; then $X \setminus A = \sigma^{-1} (Z)$. By \equref{Tg}, $R_g (u, \, u') \neq (u, \, u')$ if and only if $(u, \, u') \in Y \times Z$, and in this case $R_g (u, \, u') \in Z \times \{u\}$. Since $\sigma \times \sigma$ is injective, $R (x, \, x') \neq (x, \, x')$ if and only if $R_g (\sigma (x), \, \sigma (x')) \neq (\sigma (x), \, \sigma (x'))$, i.e. if and only if $\sigma (x) \in Y$ and $\sigma (x') \in Z$, i.e. if and only if $x \in A$ and $x' \notin A$. In this case,
\begin{equation*}
R (x, \, x') = (\sigma^{-1} \times \sigma^{-1}) \bigl( R_g (\sigma (x), \, \sigma (x')) \bigr) \in \sigma^{-1} (Z) \times \{x\} = (X \setminus A) \times \{x\}.
\end{equation*}
Thus $A$ satisfies \equref{recogTg}.
\end{proof}

\begin{remark}\relabel{preorder}
Recall that a preorder on a set $X$ is a reflexive and transitive relation $\preccurlyeq$ on $X$; it is called total if any two elements are comparable, i.e. $x \preccurlyeq x'$ or $x' \preccurlyeq x$, for all $x$, $x' \in X$. A preorder $\preccurlyeq$ defines an equivalence relation on $X$, by $x \sim x'$ if and only if $x \preccurlyeq x'$ and $x' \preccurlyeq x$; its equivalence classes are called the \emph{classes} of $\preccurlyeq$,
and $\preccurlyeq$ induces a partial order on the set of classes, which is total if $\preccurlyeq$ is total: for two distinct classes $C$ and $C'$, we say that $C$ is \emph{below} $C'$ if $c \preccurlyeq c'$ and $c' \not\preccurlyeq c$, for all $c \in C$ and $c' \in C'$.

With this terminology, \prref{recogTg} can be reformulated as follows: $(X, \, R)$ is isomorphic to a solution of the form $R_g$ if and only if the relation defined by
\begin{equation*}
x \preccurlyeq x' \qquad \text{if and only if} \qquad R (x, \, x') = (x, \, x')
\end{equation*}
is a total preorder on $X$ with at most two classes and, whenever $R (x, \, x') \neq (x, \, x')$, we have $R (x, \, x') \in [x'] \times \{x\}$, where $[x']$ is the class of $x'$. Indeed, if $A$ satisfies \equref{recogTg}, then $x \preccurlyeq x'$ if and only if $x \notin A$ or $x' \in A$; this relation is a total preorder, whose classes are those of the sets $X \setminus A$ and $A$ which are non-empty, with $X \setminus A$ below $A$, and the second part of \equref{recogTg} says that $R (x, \, x') \in [x'] \times \{x\}$. Conversely, if $\preccurlyeq$ is a total preorder with at most two classes, then either $\preccurlyeq$ is the relation $X \times X$, in which case $R = {\rm Id}_{X \times X}$ and $A := \emptyset$ satisfies \equref{recogTg}, or $\preccurlyeq$ has two classes $C_1$ and $C_2$, with $C_1$ below $C_2$, in which case $A := C_2$ satisfies \equref{recogTg}. We emphasize that neither the braid equation nor the idempotency of $R$ is assumed in \prref{recogTg}: both follow from the condition \equref{recogTg}.
\end{remark}

Using \thref{Tg}(b) we can count the solutions $R_g$ on a finite set, up to isomorphism.

\begin{corollary}\colabel{countTg}
Let $n \geq 1$. The number of isomorphism classes of solutions of the braid equation of the form $(Y \sqcup Z, \, R_g)$, with $|Y \sqcup Z| = n$, is
\begin{equation}\eqlabel{countTg}
1 + \sum_{k = 1}^{n-1} \widetilde{o} (n - k, \, k),
\end{equation}
where $\widetilde{o} (m, \, k)$ is the number of orbits of the symmetric group $S_m$, acting by conjugation, on the set of multisets with $k$ elements of maps $\{1, \ldots, m\} \to \{1, \ldots, m\}$.
\end{corollary}

\begin{proof}
If $Y = \emptyset$ or $Z = \emptyset$, then $R_g$ is the identity map; this gives the term $1$ in \equref{countTg}. Assume now that $Y$ and $Z$ are non-empty; then $R_g$ is not the identity map, since it moves all the pairs of $Y \times Z$. By the proof of \thref{Tg}(b), any isomorphism $(Y \sqcup Z, \, R_g) \to (Y' \sqcup Z', \, R_{g'})$ maps $Y$ onto $Y'$; hence solutions with different values of $k := |Y|$ are not isomorphic. Let $1 \leq k \leq n - 1$. Up to isomorphism, we can assume that $Z = \{1, \ldots, n-k\}$, for all the triples $(Y, \, Z, \, g)$ with $|Y| = k$. For such a triple, let $M_g$ be the multiset of values of $g$, i.e. the multiset with the underlying set $\{g_y \, | \, y \in Y\}$, in which $h$ has the multiplicity $|\{y \in Y \, | \, g_y = h\}|$. Then $M_g$ is a multiset with $k$ elements of maps $Z \to Z$, and any such multiset is of the form $M_g$. The condition \equref{isoTg} can be written as $g'_{\alpha (y)} = \beta \circ g_y \circ \beta^{-1}$, for all $y \in Y$; since $\alpha$ is an arbitrary bijection, for a given $\beta \in S_{n-k}$ there exists a bijection $\alpha \colon Y \to Y'$ satisfying this condition if and only if $\beta \cdot M_g = M_{g'}$. Thus, by \thref{Tg}(b), $R_g \cong R_{g'}$ if and only if $M_g$ and $M_{g'}$ belong to the same orbit of $S_{n-k}$, and the number of isomorphism classes of solutions with $|Y| = k$ is $\widetilde{o} (n - k, \, k)$. This proves the formula \equref{countTg}.
\end{proof}

\begin{remark}\relabel{countTg}
The number $\widetilde{o} (m, \, k)$ is the analogue of the number $o (m, \, k)$ used in \cite[Corollary 2.13]{MilBIP} for counting the bipartite solutions, where only multisets of permutations are considered; in particular, $o (m, \, k) \leq \widetilde{o} (m, \, k)$.

For $k = 1$, $\widetilde{o} (m, \, 1) = d (m)$ is the Davis number \cite[Theorem 6]{davis} of a set with $m$ elements, i.e. the number of conjugacy classes of maps $Z \to Z$ on a set $Z$ with $m$ elements, or, equivalently, the number of functional graphs with $m$ vertices, up to isomorphism \cite[Remark 1.4]{acm1}; for $m = 1, \ldots, 5$ we have $d (m) = 1$, $3$, $7$, $19$, $47$ (sequence A001372 of \cite{OEIS}).

The numbers $\widetilde{o} (m, \, k)$ can be computed by the Cauchy--Frobenius Lemma \cite[Theorem 3.22]{Rotman}. For $\sigma \in S_m$, let $c_j (\sigma)$ be the number of cycles of length $j$ of $\sigma$, for $j \geq 1$, and let
\begin{equation*}
N (\sigma) := \prod_{j \geq 1} \Bigl( \sum_{i \mid j} i \, c_i (\sigma) \Bigr)^{c_j (\sigma)}
\end{equation*}
be the number of maps $f \colon \{1, \ldots, m\} \to \{1, \ldots, m\}$ which commute with $\sigma$ \cite[Theorem 2.18(3)]{Daug}: such a map is determined on each cycle of length $j$ of $\sigma$ by the image of one of its points, which can be any point lying on a cycle of length dividing $j$. Then $\widetilde{o} (m, \, k)$ is the coefficient of $t^k$ in the power series
\begin{equation}\eqlabel{frobenius}
\frac{1}{m!} \sum_{\sigma \in S_m} \exp \Bigl( \sum_{r \geq 1} \frac{N (\sigma^r)}{r} \, t^r \Bigr).
\end{equation}
Indeed, by the Cauchy--Frobenius Lemma, $\widetilde{o} (m, \, k) = \frac{1}{m!} \sum_{\sigma \in S_m} F_k (\sigma)$, where $F_k (\sigma)$ is the number of multisets with $k$ elements of maps which are fixed by $\sigma$, i.e. whose multiplicity function (extended by $0$) is constant on the orbits of $\langle \sigma \rangle$ acting by conjugation on the set of all maps. If these orbits have the lengths $\ell_1, \ldots, \ell_p$, then
\begin{equation*}
\sum_{k \geq 0} F_k (\sigma) \, t^k = \prod_{s = 1}^{p} \frac{1}{1 - t^{\ell_s}} = \exp \Bigl( \sum_{r \geq 1} \frac{t^r}{r} \sum_{s: \, \ell_s \mid r} \ell_s \Bigr) = \exp \Bigl( \sum_{r \geq 1} \frac{N (\sigma^r)}{r} \, t^r \Bigr),
\end{equation*}
since $\sum_{s: \, \ell_s \mid r} \ell_s$ is the number of maps fixed by the conjugation with $\sigma^r$, i.e. $N (\sigma^r)$. In other words, \equref{frobenius} is a particular case of P\'olya's enumeration theorem \cite[Chapter 2]{HP73}. For $k = 1$, \equref{frobenius} gives de Bruijn's formula $d (m) = \frac{1}{m!} \sum_{\sigma \in S_m} N (\sigma)$ \cite[(1.1)]{acm1}, written as a sum over permutations.

A computer program based on \equref{frobenius} shows that, for sets with $n = 1, \ldots, 10$ elements, the formula \equref{countTg} gives $1$, $2$, $5$, $16$, $108$, $2183$, $165654$, $51696309$, $57169960202$, $313513723362924$ isomorphism classes of solutions of the form $R_g$.
\end{remark}

%%%%%%%%%%%%%%%%%%%%%%%%%%%%%%%%%%%%%%%%%%%%%%%%%%%%%%%%%%%%%%%%%%%%%%%%%%%%%
\section{Centralizers of finitely many maps and the braid equation}\selabel{centr}

Recall from \cite[Definition 1.1]{acm1} that the \emph{centralizer} of a family $(f_i)_{i\in I}$ of self-maps $f_i \colon X \to X$ is the subgroup of ${\rm Sym} (X)$ defined by
\begin{equation*}
{\rm Aut} (f_i, \, i\in I) := \bigl\{ \sigma \in {\rm Sym} (X) \, \bigl| \, \sigma \circ f_i = f_i \circ \sigma, \,\, \text{for all} \,\, i \in I \bigr\}.
\end{equation*}
It was proved in \cite[Theorem 2.19]{acm1} that the centralizer ${\rm Aut} (f)$  of a single map $f \colon X \to X$ is the automorphism group of a left non-degenerate solution
of the braid equation. On the other hand, the groups which are isomorphic to the centralizer ${\rm Aut} (f)$ of a single map $f \colon X \to X$ are classified in \cite[Section 4]{acm1} and summarized in Theorem F of the Introduction of \cite{acm1}; in particular, the abelian groups of the form ${\rm Aut} (f)$ are exactly the direct products of cyclic groups \cite[Corollary 4.22]{acm1}. Consequently, any finitely generated abelian group is isomorphic to ${\rm Aut} (f)$, for some map $f$, while the additive group $(\bQ, +)$ of rational numbers is not.

Perhaps surprisingly, any group $G$ is isomorphic to the centralizer of a family of two maps \cite[Theorem 4.2(1)]{acm1}, i.e. $G \cong {\rm Aut} (f_1, \, f_2)$, for some set $X$ and two
maps $f_1$, $f_2 \colon X \to X$. In this section we show, using the solutions $R_g$ of \seref{Tg}, that the centralizer of any finite family of maps $f_1, \ldots, f_k \colon X \to X$ is the automorphism group of an idempotent solution of the braid equation. In particular, we obtain another proof of the fact that any group is isomorphic to the automorphism group of a solution of the braid equation.

Let $X$ be a set and $f = (f_1, \ldots, f_k)$ a finite family of maps $f_i \colon X \to X$, $k \geq 1$. We consider the sets
\begin{equation*}
Y_f := \{y_0, \, y_1, \ldots, y_k\}, \qquad W := \{w_0, \, w_1, \ldots, w_{k+1}\}, \qquad Z_f := X \sqcup W,
\end{equation*}
where $Y_f$ and $W$ consist of $k+1$, respectively $k+2$, new elements, and the maps $g_{y_i} \colon Z_f \to Z_f$, $i = 0, \ldots, k$, defined for $x\in X$ and $j = 1, \ldots, k+1$ by
\begin{alignat*}{3}
g_{y_0} (x) &:= w_0, &\qquad g_{y_0} (w_0) &:= w_0, &\qquad g_{y_0} (w_j) &:= w_{j-1},\\
g_{y_i} (x) &:= f_i (x), & g_{y_i} (w_i) &:= w_0, & g_{y_i} (w_l) &:= w_l \,\,\, \text{for} \,\, l \neq i, \qquad \text{if} \,\, 1 \leq i \leq k.
\end{alignat*}

We denote by $R_f := R_g$ the idempotent solution of the braid equation on $X_f := Y_f \sqcup Z_f = Y_f \sqcup X \sqcup W$ given by \thref{Tg} for the family of maps $g = (g_{y_i})_{0 \leq i \leq k}$.

\begin{theorem}\thlabel{centr}
Let $f_1, \ldots, f_k \colon X \to X$ be a finite family of maps, $k \geq 1$. Then the map
\begin{equation*}
{\rm Aut} (f_1, \ldots, f_k) \to {\rm Aut} (X_f, \, R_f), \qquad \gamma \mapsto {\rm Id}_{Y_f} \sqcup \gamma \sqcup {\rm Id}_W,
\end{equation*}
is an isomorphism of groups. In particular, the centralizer of $f_1, \ldots, f_k$ is the automorphism group of an idempotent solution of the braid equation, and of the YB-solution $\tau_{X_f} \circ R_f$, on a set with $|X| + 2k + 3$ elements.
\end{theorem}

\begin{proof}
Since $Y_f$ and $Z_f$ are non-empty, \thref{Tg}(c) shows that the automorphisms of $(X_f, \, R_f)$ are exactly the bijections $\alpha \sqcup \beta$, where $\alpha \in {\rm Sym} (Y_f)$ and $\beta \in {\rm Sym} (Z_f)$ satisfy
\begin{equation}\eqlabel{conjg}
\beta \circ g_{y_i} = g_{\alpha (y_i)} \circ \beta, \qquad \text{for all} \,\, i = 0, \ldots, k.
\end{equation}

First, let $\gamma \in {\rm Aut} (f_1, \ldots, f_k)$ and $\beta := \gamma \sqcup {\rm Id}_W \in {\rm Sym} (Z_f)$. We show that the pair $({\rm Id}_{Y_f}, \, \beta)$ satisfies \equref{conjg}, i.e. $\beta \circ g_{y_i} = g_{y_i} \circ \beta$, for all $i$. On $W$ both sides are equal to $g_{y_i}$, since $\beta$ is the identity on $W$ and $g_{y_i} (W) \subseteq W$. For $x \in X$ we have
\begin{equation*}
\beta \bigl( g_{y_0} (x) \bigr) = \beta (w_0) = w_0 = g_{y_0} \bigl( \gamma (x) \bigr) = g_{y_0} \bigl( \beta (x) \bigr),
\end{equation*}
and, for $1 \leq i \leq k$,
\begin{equation*}
\beta \bigl( g_{y_i} (x) \bigr) = \gamma \bigl( f_i (x) \bigr) = f_i \bigl( \gamma (x) \bigr) = g_{y_i} \bigl( \beta (x) \bigr).
\end{equation*}
Hence ${\rm Id}_{Y_f} \sqcup \gamma \sqcup {\rm Id}_W$ is an automorphism of $(X_f, \, R_f)$. Clearly, the map $\gamma \mapsto {\rm Id}_{Y_f} \sqcup \gamma \sqcup {\rm Id}_W$ is an injective morphism of groups.

Conversely, let $\sigma = \alpha \sqcup \beta$ be an automorphism of $(X_f, \, R_f)$, where $\alpha$ and $\beta$ satisfy \equref{conjg}. Since $\beta$ is a bijection, \equref{conjg} can be written as $g_{\alpha (y_i)} = \beta \circ g_{y_i} \circ \beta^{-1}$; in particular, $(g_{\alpha (y_i)})^n = \beta \circ (g_{y_i})^n \circ \beta^{-1}$, for all $n \geq 1$. We show, in three steps, that $\alpha = {\rm Id}_{Y_f}$, that $\beta$ is the identity on $W$, and that the restriction of $\beta$ to $X$ belongs to ${\rm Aut} (f_1, \ldots, f_k)$.

Step 1. We first show that $\alpha (y_0) = y_0$. The map $(g_{y_0})^{k+1}$ is constant, equal to $w_0$: indeed, $g_{y_0}$ maps every element of $X \cup \{w_0\}$ to $w_0$, and $w_j$ to $w_{j-1}$, for $1 \leq j \leq k+1$. On the other hand, if $1 \leq i \leq k$, then $g_{y_i}$ fixes the two distinct elements $w_0$ and $w_{k+1}$; hence no power of $g_{y_i}$ is constant. As a conjugate of a constant map by a bijection is again constant, $(g_{\alpha (y_0)})^{k+1} = \beta \circ (g_{y_0})^{k+1} \circ \beta^{-1}$ is constant. Therefore $\alpha (y_0) = y_0$, and \equref{conjg} for $i = 0$ shows that $\beta$ commutes with $g_{y_0}$:
\begin{equation*}
\beta \circ g_{y_0} = g_{y_0} \circ \beta.
\end{equation*}
This will be used in Step 2. Moreover, since $\alpha$ is a bijection, $\alpha$ maps the set $\{y_1, \ldots, y_k\}$ onto itself; this will be used in Step 3.

Step 2. Next, we show that $\beta$ is the identity on $W$; in particular, $\beta (X) = X$. We use only the fact that $\beta$ commutes with $g_{y_0}$. Since $w_0$ is the unique fixed point of $g_{y_0}$, $\beta (w_0)$ is also a fixed point of $g_{y_0}$, i.e. $\beta (w_0) = w_0$. Next, the image of $g_{y_0}$ is the set $\{w_0, \, w_1, \ldots, w_k\}$, and it is invariant under $\beta$, since $\beta (g_{y_0} (Z_f)) = g_{y_0} (\beta (Z_f)) = g_{y_0} (Z_f)$. The elements of this image which are mapped to $w_0$ by $g_{y_0}$ are $w_0$ and $w_1$. Now $\beta (w_1)$ belongs to the image of $g_{y_0}$, $\beta (w_1) \neq \beta (w_0) = w_0$, and
\begin{equation*}
g_{y_0} \bigl( \beta (w_1) \bigr) = \beta \bigl( g_{y_0} (w_1) \bigr) = \beta (w_0) = w_0;
\end{equation*}
hence $\beta (w_1) = w_1$. Finally, for $1 \leq j \leq k$, the element $w_{j+1}$ is the unique $z \in Z_f$ such that $g_{y_0} (z) = w_j$. If $\beta (w_j) = w_j$, then $g_{y_0} (\beta (w_{j+1})) = \beta (g_{y_0} (w_{j+1})) = \beta (w_j) = w_j$, hence $\beta (w_{j+1}) = w_{j+1}$. By induction, $\beta$ fixes $W$ pointwise and, as $\beta$ is a bijection of $Z_f = X \sqcup W$, we obtain $\beta (X) = X$.

Step 3. Finally, we show that $\alpha = {\rm Id}_{Y_f}$ and that $\gamma := \beta|_X$ belongs to ${\rm Aut} (f_1, \ldots, f_k)$. Let $1 \leq i \leq k$. By Step 1, $\alpha (y_i) = y_l$, for some $1 \leq l \leq k$; note that $l = 0$ is excluded, which is essential, since $g_{y_0} (w_1) = w_0$. Evaluating \equref{conjg} at $w_i$ and using Step 2, we obtain
\begin{equation*}
w_0 = \beta (w_0) = \beta \bigl( g_{y_i} (w_i) \bigr) = g_{y_l} \bigl( \beta (w_i) \bigr) = g_{y_l} (w_i).
\end{equation*}
Since $g_{y_l} (w_i) = w_i \neq w_0$ for $l \neq i$, we get $l = i$. Thus $\alpha = {\rm Id}_{Y_f}$, and \equref{conjg} becomes $\beta \circ g_{y_i} = g_{y_i} \circ \beta$, for all $i$. By Step 2, $\gamma := \beta|_X$ is a bijection of $X$ and, for $x \in X$ and $1 \leq i \leq k$,
\begin{equation*}
\gamma \bigl( f_i (x) \bigr) = \beta \bigl( g_{y_i} (x) \bigr) = g_{y_i} \bigl( \beta (x) \bigr) = f_i \bigl( \gamma (x) \bigr).
\end{equation*}
Hence $\gamma \in {\rm Aut} (f_1, \ldots, f_k)$ and $\sigma = {\rm Id}_{Y_f} \sqcup \gamma \sqcup {\rm Id}_W$. This proves that the map of the statement is an isomorphism of groups.

Finally, $|X_f| = |Y_f| + |X| + |W| = |X| + 2k + 3$, and the statement about the YB-solution $\tau_{X_f} \circ R_f$ follows from \equref{auttau}.
\end{proof}

\begin{corollary}\colabel{centr}
Any group $G$ is isomorphic to the automorphism group of an idempotent solution of the braid equation of the form $R_f$, for a family $f = (f_1, \, f_2)$ of two maps.
\end{corollary}

\begin{proof}
By \cite[Theorem 4.2(1)]{acm1}, $G \cong {\rm Aut} (f_1, \, f_2)$, for two maps $f_1$, $f_2 \colon X \to X$ on some set $X$. Now apply \thref{centr} for $k = 2$.
\end{proof}

\begin{example}\exlabel{centr1}
Let $f \colon X \to X$ be a map, i.e. $k = 1$ in \thref{centr}. Then $Y_f = \{y_0, \, y_1\}$, $W = \{w_0, \, w_1, \, w_2\}$ and $X_f = \{y_0, \, y_1\} \sqcup X \sqcup \{w_0, \, w_1, \, w_2\}$ is a set with $|X| + 5$ elements. The idempotent solution $R_f \colon X_f \times X_f \to X_f \times X_f$ of the braid equation given by \equref{Tg} is defined for all $x \in X$ by
\begin{alignat*}{4}
R_f (y_0, \, x) &= (w_0, \, y_0), &\quad R_f (y_0, \, w_0) &= (w_0, \, y_0), &\quad R_f (y_0, \, w_1) &= (w_0, \, y_0), &\quad R_f (y_0, \, w_2) &= (w_1, \, y_0),\\
R_f (y_1, \, x) &= \bigl( f (x), \, y_1 \bigr), & R_f (y_1, \, w_0) &= (w_0, \, y_1), & R_f (y_1, \, w_1) &= (w_0, \, y_1), & R_f (y_1, \, w_2) &= (w_2, \, y_1),
\end{alignat*}
and $R_f (a, \, b) = (a, \, b)$, for all the other pairs $(a, \, b) \in X_f \times X_f$. By \thref{centr},
\begin{equation*}
{\rm Aut} (X_f, \, R_f) = \bigl\{ {\rm Id}_{Y_f} \sqcup \gamma \sqcup {\rm Id}_W \, \bigl| \, \gamma \in {\rm Aut} (f) \bigr\} \cong {\rm Aut} (f).
\end{equation*}
For instance, let $n \geq 2$, $X = \bZ_n$ and $f (x) = x + 1$, for all $x \in \bZ_n$. A permutation $\gamma$ of $\bZ_n$ commutes with $f$ if and only if $\gamma (x) = \gamma (f^x (0)) = f^x (\gamma (0)) = x + \gamma (0)$, for all $x \in \bZ_n$; hence ${\rm Aut} (f) = \{ x \mapsto x + v \, | \, v \in \bZ_n \} \cong C_n$. Thus the cyclic group $C_n$ is the automorphism group of an idempotent solution of the braid equation on a set with $n + 5$ elements.

However, \thref{centr} does not give the smallest set on which ${\rm Aut} (f)$ is realized as an automorphism group: by \cite[Theorem 2.19]{acm1}, ${\rm Aut} (f)$ is the automorphism group of the left non-degenerate FS-solution (in particular, a solution of the braid equation) $R \colon X \times X \to X \times X$, $R (x, \, y) = (y, \, f (x))$, which lives on the set $X$ itself, while $R_f$ lives on a set with $|X| + 5$ elements.
\end{example}

\begin{remark}\relabel{smallest}
The set $X$ is not always the smallest one on which ${\rm Aut} (f)$ can be realized as the automorphism group of a solution of the braid equation. For instance, if $f$ is the constant map with value $c$ on a set $X$ with $m \geq 2$ elements, then ${\rm Aut} (f) = \{ \sigma \in {\rm Sym} (X) \, | \, \sigma (c) = c \} \cong S_{m-1}$, which is also the automorphism group of the identity map on a set with $m - 1$ elements. On the other hand, since ${\rm Aut} (X, \, R)$ is a subgroup of ${\rm Sym} (X)$, a finite group $G$ is isomorphic to the automorphism group of a solution $R$ on a set $X$ only if $|X| \geq \mu (G)$, where $\mu (G)$ is the minimal degree of a faithful permutation representation of $G$. For instance, $\mu (C_n) = n$ if $n$ is a prime power; in this case the solution $R (x, \, y) = (y, \, x + 1)$ on $\bZ_n$ lives on the smallest possible set.

In general, finding, for a given finite group $G$, the smallest set on which $G$ is the automorphism group of a solution of the braid equation is an open problem (see \cite[Section 3]{MilAut}).
\end{remark}

%%%%%%%%%%%%%%%%%%%%%%%%%%%%%%%%%%%%%%%%%%%%%%%%%%%%%%%%%%%%%%%%%%%%%%%%%%%%%

\subsection*{Acknowledgements}
We used Anthropic, Claude Opus~5.5, during the development and drafting of the paper to explore possible formulations, proof strategies, proofreading, correcting typos, and making improvements. The Python program used in \reref{countTg} to compute the numbers $\widetilde{o} (m, \, k)$ was written by Claude Opus~5.5; the author assumes full responsibility for it, and the program is available from the author upon request. All statements, proofs and computational results were independently checked by the author, who takes sole responsibility for them.

\newcommand{\Addresses}{{%
  \bigskip
  \footnotesize
  \textsc{Faculty of Mathematics and Computer Science, University of Bucharest}
  \par\nopagebreak
  \textsc{Str. Academiei 14, RO-010014 Bucharest 1, Romania}
  \par\nopagebreak
  \textit{E-mail address}: \texttt{gigel.militaru@gmail.com}
}}

\Addresses

\end{document}